\documentclass[11pt, oneside, a4paper]{article}
\usepackage{amssymb,amscd}
\usepackage{amsmath}
\usepackage{xypic}
\newtheorem{theorem}{Theorem}[section]

\newtheorem{lemma}[theorem]{Lemma}
\newtheorem{cor}[theorem]{Corollary}
\newtheorem{prop}[theorem]{Proposition}

\newtheorem{defn}[theorem]{Definition}

\def\im{{\rm im}}

\newcommand{\Om}{\Omega}

\newcommand{\OmK}{[\Omega]^k}

\newcommand{\cvd}{\hspace*{\fill}
    {\rm \hbox{\vrule height 0.2 cm width 0.2cm}}}
\newcommand{\N}{{\mathbb N}}
\newcommand{\Fi}{{\mathbb F}}

\newcommand{\Z}{{\mathbb Z}}

 \DeclareMathOperator{\Sym}{Sym}

\DeclareMathOperator{\Ker}{ker}

\title{Second cohomology groups and finite covers}

\author{David M. Evans\\
 School of Mathematics\\
University of East Anglia, Norwich, NR47TJ, UK\\
\vspace{0.5cm}
{\tt{d.evans@uea.ac.uk}}\\
Elisabetta Pastori\\
Dipartimento di Matematica\\
Universit\`a degli Studi di Torino, Via Carlo Alberto, 10 10123\\
Torino, Italy\\
{\tt{elisabetta.pastori@unito.it}}}

\begin{document}
\date{}

\maketitle
\begin{abstract}
For $\Omega$ an infinite set,  $k \geq 2$ and $W$ the set of $k$-sets from $\Omega$, there is a natural closed permutation group $\Gamma_k$ which is a non-split extension $0 \to \Z_2^W \to \Gamma_k \to \Sym(\Omega) \to 1$. We classify the closed subgroups of $\Gamma_k$ which project onto $\Sym(\Omega)$. The question arises in model theory as a problem about finite covers, but here we formulate and solve it in algebraic terms.
\end{abstract}

\section{Introduction and background on finite covers}

The problem of understanding the \emph{finite covers}  of a structure arises in model theory (for example, see \cite{AZ2}, \cite{HPi}), but  it also has a natural formulation in purely algebraic terms as an extension problem in the category of permutation groups and we adopt this approach here. We begin by reviewing some definitions and basic facts about infinite permutation groups and  finite covers. Further details can be found in \cite{EMI} and \cite{EHew}. 

If $C$ is any set then the full symmetric group $\Sym(C)$  on $C$ can be considered as a topological group by giving it the topology whose  open sets are arbitrary unions of cosets of pointwise stabilizers of finite subsets of $C$.  A subgroup $\Gamma\leq \Sym(C)$ is closed in $\Sym(C)$ if and only if any element of $\Sym(C)$ which stabilizes each $\Gamma$-orbit on $C^n$, for all $n\in \N$, is in $\Gamma$. As  is well known,  closed subgroups in this topology are precisely automorphism groups of first-order structures on $C$.  Hence  we say that a \emph{permutation structure} is a pair $\langle C; \Gamma\rangle$ where $C$ is a  set and $\Gamma$ is a closed subgroup of $\Sym(C)$.

Suppose $\langle C;\Gamma \rangle$ is a permutation structure and $\rho$ is a $\Gamma$-invariant equivalence relation on $C$. Let $W=C/\rho$, the set of equivalence classes. The action of $\Gamma$ on $W$ gives a continuous homomorphism $\mu:\Gamma\rightarrow \Sym(W)$. If all the $\rho$-classes are finite then $G=\mu(\Gamma)$ is a closed subgroup of $\Sym(W)$, $\mu$ is an open map (indeed a proper map) and $\langle W; G\rangle$ is a permutation structure. We say that  $\langle C;\Gamma \rangle$ is a \emph{finite cover} of  $\langle W; G\rangle$ with \emph{kernel} $\Ker\mu$. The latter is a profinite group and we have an exact sequence of topological groups
\begin{equation}\label{exact}
1\rightarrow K\rightarrow \Gamma \stackrel{\mu}{\rightarrow} G\rightarrow 1.
\end{equation}

We say that a  subgroup $\Sigma$ of $\Gamma$ is \textit{full} (with respect to $\mu$) if it is closed and $\mu(\Sigma) = G$, or equivalently, if $\Gamma = K\Sigma$. Usually $\mu$ will be clear from the context and in this case we will just say that $\Sigma$ is full in $\Gamma$. We say that $\Sigma$ is \textit{minimally full} if it is full and has no proper full subgroup. It can be shown that any full subgroup of $\Gamma$ contains a minimally full subgroup (Lemma 1.5(i) of \cite{CKK}). The main result of this paper is an explicit determination of the minimally full subgroups of some very natural finite covers associated with the infinite symmetric group. To describe these and to motivate the problem a little further,  it will be useful to introduce the following terminology.

A finite cover $\langle C; \Gamma \rangle$ of $\langle W; G\rangle$ with kernel $K$ has the following data associated to it. Let $w\in W$ and $C(w)$ denote the $\rho$-class in $C$ labelled by $w$. For $w\in W$, the \emph{fibre group} $F(w)$ at  $w$ is the finite permutation group induced on $C(w)$ by the setwise stabilizer in $\Gamma$ of $C(w)$. The  \emph{binding group} $B(w)$ at $w$ is the finite permutation group induced on $C(w)$ by the kernel $K$. It is clear that $B(w)\unlhd F(w)$.  Let $G_w$ be the stabilizer of $w$ in $G$. There is a homomorphism $\chi_w: G_w\rightarrow F(w)/B(w)$,
defined by $\chi_w(g)=(h|C(w))B(w),$ where $g\in G_w$
and $h\in \Gamma$ is a permutation which extends $g$.
This homomorphism is well defined, continuous and surjective
(\cite{EMI}, Lemma 2.1.1). We shall refer to it as the
{\emph{canonical
homomorphism}} of the cover.  Clearly, if  $G$ acts transitively on $W$ then all the fibre groups are
isomorphic as permutation groups, as are the
binding groups. 

 Conversely, given a permutation structure $\langle W; G\rangle$ with $G$ acting transitively on $W$, some $w \in W$, finite permutation groups $B\unlhd F$ and a continuous epimomorphism $\chi : G_w\rightarrow F/B$, there is a unique finite cover $\langle C_0; \Gamma_0\rangle$ of $\langle W; G\rangle$ with these data and kernel $K_0$ equal to the full direct product $\prod_{w\in W}B$ (\cite{EMI}, Lemma 2.1.2). We refer to this as the \textit{free} finite cover with the given data.  Any finite cover $\langle C; \Gamma\rangle$ of $\langle W, G\rangle$ with the same data can be regarded as having $C = C_0$ and $\Gamma \leq \Gamma_0$. 
 
So the problem of determining the finite covers of a transitive permutation structure $\langle W; G\rangle$ with given fibre and binding groups and canonical homomorphisms is equivalent to the problem of determining the full subgroups $\Gamma$ of the corresponding free cover $\langle C_0; \Gamma_0\rangle$. Moreover we can write $\Gamma = K \Gamma_1$ where $K = \Gamma \cap \Gamma_0$ is the kernel of the cover and $\Gamma_1$ is some minimally full subgroup. Thus the problem can be reduced to determining the possible kernels $K \leq K_0$ and the  minimally full subgroups of $\Gamma_0$. In the case where the binding groups are abelian, $K_0$ is a continuous $G$-modules and the possible kernels are just the closed $G$-submodules. 

Previous work on problems of this kind have been in cases where the fibre groups split over the binding groups. Our main results here are concerned with the following finite covers, which are probably the simplest examples where this is not the case. Denote by $\Z_n$ the cyclic group of order $n$.

\begin{defn}\rm \label{ckdef}
Let $\Omega$ be an infinite set and $k \geq 2$ a natural number. Then $G = \Sym(\Omega)$ acts naturally on  $[\Om]^k$, the set of $k$-subsets from $\Omega$ and is a closed permutation group on this set. Note that if $w \in [\Om]^k$, then $G_w = \Sym(w) \times \Sym(\Om \setminus w)$ so there is a homomorphism $\chi_w : G_w \to \Z_2$ given by taking $\chi_w(g)$ to be the sign of $g$ restricted to $w$. Let $\langle C_k; \Gamma_k\rangle$ be the free finite cover of $\langle [\Om]^k; G\rangle$ with fibre groups $\Z_4$ acting regularly, binding group $\Z_2$, and canonical homomorphisms $\chi_w$.
\end{defn}

So with this notation, we have an exact sequence:
\[ 0 \to \Z_2^{[\Om]^k} \to \Gamma_k \to \Sym(\Om) \to 1. \]
The paper \cite{Gr} gives a method for determining the closed $G$-submodules of $\Z_2^{[\Om]^k}$: cf. Section \ref{CG} here. Thus, to determine all finite covers of $\langle [\Om]^k; G\rangle$ with the given data it remains to determine the minimally full subgroups of $\Gamma_k$. We summarise our results on this in the following.

\begin{theorem}
\begin{enumerate}
\item[(i)] The group $\Gamma_2$ has no proper full subgroup.
\item[(ii)] For $2 \leq \ell$ there exist continuous homomorphisms $\gamma_{2,\ell} : \Gamma_2 \to \Gamma_\ell$ which extend the natural $G$-homomorphisms $\alpha_{2,\ell} : \Z_2^{[\Om]^2} \to \Z_2^{[\Om]^\ell}$. The $\Gamma_\ell$-conjugates of $\gamma_{2,\ell}(\Gamma_2)$ are the minimally full subgroups of $\Gamma_\ell$.
\end{enumerate}
\end{theorem}

Part (i) here is Theorem \ref{application1} and answers Question 8.8 of \cite{EMI}. Part (ii) is Theorem \ref{47}.

\section{Full subgroups and continuous cohomology}

We use cohomological methods to anaylse group extensions as in the exact sequence (1) of the Introducton. As these are topological groups, the appropriate context is that of continuous cohomology. The basic results of group cohomology (extension theory, the long exact sequence, Shapiro's lemma) hold in this context, although some care is needed in adapting the standard proofs. For  more details we refer the reader to \cite{EHew}.  Previous work on finite covers has made use of first cohomology groups $H^1_c$. One of the novelties of our approarch here is that we use second-degree cohomology $H^2_c$.

Denote by $\mathcal{PG}$ the class of all topological groups aring as permutation groups. So these are Hausdorff topological groups with a base of open neighbourhoods of the identity consisting of open subgroups. Suppose $G \in \mathcal{PG}$ is a closed permutation group and $K$ is a continuous profinite $G$-module.  We denote by  $C_c^n(G,K)$   the additive group of continuous functions $\varphi:G^n \rightarrow K$. The usual coboundary operator $\delta^n$ sends $C_c^n(G,K)$ to $ C_c^{n+1}(G,K)$, so that $(C_c^n(G,K);\delta^n)_{n\in \mathbb{N}}$ is a cochain complex. The homology of this complex, $H_c^\ast(G,K)$, is the  \emph{ continuous cohomology} of $G$ with coefficients in $K$. 
 
 Let $G\in \mathcal{PG}$ and $M$ be a profinite continuous $G$-module (written additively). By a \emph{$\mathcal{PG}$-extension} of $M$ by
$G$ we mean a short exact sequence with continuous open
homomorphisms
\begin{equation}\label{shortexact}
0\rightarrow M \rightarrow \Gamma
\rightarrow G \rightarrow 1,
\end{equation}
where  $\Gamma \in \mathcal{PG}$.

Any such $\mathcal{PG}$-extension (\ref{shortexact}) admits a continuous closed section $s:G\rightarrow \Gamma$ (\cite{EHew}) and so cohomology of low degree continuous cocycles on profinite $G$-modules retains its familiar applications: $H^1_c(G,M)$ classifies closed complements in the split extension and $H^2_c(G,M)$ classifies all $\mathcal{PG}$-extensions of $M$ by $G$ (see the proof of Theorem \ref{main theorem} below).

We use the following notation. If $f:M_1\rightarrow M_2$ is a continuous $G$-homomorphism of  profinite continuous $G$-modules, we denote by $$
f^\ast:H_c^n(G,M_1)\rightarrow H_c^n(G,M_2)
$$
the induced homomorphism in cohomology.  Note if $M_0$ is the kernel of $f$, this is the map in the long exact sequence

$$
\cdots \rightarrow H_c^{n-1}(G,M_0)\rightarrow H_c^{n-1}(G,M_1)\rightarrow H_c^{n-1}(G,M_2)\rightarrow H_c^{n}(G,M_0)\rightarrow \cdots
$$
arising from the short exact sequence $0 \to M_0 \to M_1 \to M_2 \to 0$. 

\begin{theorem}\label{main theorem}
Suppose $\langle C_0; \Gamma_0\rangle$ is a finite cover of $\langle W; G\rangle$ with abelian kernel $K_0$ and $K$ a closed $G$-submodule of $K_0$. Let $e_0$ be the element
in $H^2_c(G,K_0)$ which gives rise to $\Gamma_0$ (as an extension of
$K_0$ by $G$).

Let $$ 0\rightarrow K\xrightarrow{i} K_0\rightarrow \bar{K}
\rightarrow 0
$$
be the natural short exact sequence where $i$ is the inclusion map. Consider 
$$
\cdots \rightarrow H^1_c(G,\bar{K})\rightarrow
H^2_c(G,K)\xrightarrow{i^\ast} H^2_c(G,K_0),
$$
 part of the long exact sequence, where $i^\ast$ is the induced map in
cohomology. Then there exists a full subgroup $\Gamma \leq \Gamma_0$ with
 $K = \Gamma\cap K_0$ if and only if there exists an element $e\in
H^2_c(G,K)$ such that $i^\ast(e)=e_0$.\end{theorem}
\emph{Proof.} 
We adapt some standard facts about group extensions, omitting the details of the proofs. Let $0\rightarrow N\rightarrow E\rightarrow L\rightarrow 1$ be a $\mathcal{PG}$-extension and $s: L\rightarrow  E$ a continuous closed section. Given $g_1,g_2\in L$, since $s(g_1)+s(g_2)$
and $s(g_1g_2)$ belong to the same coset of $N$ in $E$,
there exists  an element $f(g_1,g_2)\in N$ such that
\begin{equation}\label{factor set}
s(g_1)+s(g_2)=f(g_1,g_2)+s(g_1g_2).
\end{equation} The map $f:L^2\rightarrow N$  is a continuous 2-cocycle. We can  then define a multiplication on the set $L\times N$ by \begin{equation}\label{moltiplication}
(g_1,n_1)(g_2,n_2)=(g_1g_2,f(g_1,g_2)+g_2 n_1+n_2).
\end{equation}
Then $L\times N$ with  multiplication (\ref{moltiplication}) becomes a topological group and $0\rightarrow N \rightarrow L\times N
\rightarrow L \rightarrow 1$ is a $\mathcal{PG}$-extension equivalent to the original one.  All
$\mathcal{PG}$-extensions of $N$ by $L$  are obtained in this way up to isomorphism, and
two $\mathcal{PG}$-extensions are isomorphic if and only if their continuous
2-cocycles differ by a continuous 2-coboundary. Also, different choices of section give cohomologous continuous 2-cocycles.

Suppose now that $\Gamma$ is a full subgroup of $\Gamma_0$ with $\Gamma\cap K_0 = K$.  So $\Gamma$ is
a $\mathcal{PG}$-extension of $K$ by $G$. Let $s:  G  \rightarrow  \Gamma$, given by $s(a)=g_a$, be a continuous closed section associated to the $\mathcal{PG}$-extension $ \Gamma$  of $K$ by $G$ and let $m_K:G\times G\rightarrow K$ be the corresponding (given by (\ref{factor set})) continuous
2-cocycle. We consider the map
$$
\begin{array}{cccccc}
  s_0: & G & \overset{s}{\rightarrow} & \Gamma & \overset{i}{\hookrightarrow} & \Gamma_0 \\
   & a & \mapsto & g_a & \mapsto & g_a\\
\end{array}
$$
which is continuous and closed. It is easy to check that $s_0:G\rightarrow \Gamma_0$ is a  section. We also denote the resulting 2-cocycle (which goes from $G\times G$ to $K_0$) by $m_K$ since  on each element of $G\times G$ it takes the same value as the 2-cocycle  $m_K:G\times G\rightarrow K$.
It follows that $e_0=m_{K}+ B^2_c(G,K_0)$ and, if $e$ is the element in $H^2_c(G,K)$  which $\Gamma$ comes from, then
$e=m_K+B^2_c(G,K)$. Hence, $i^\ast(e)=e_0$.\\

Conversely, suppose $e=m_K+B^2_c(G,K)\in
H^2_c(G,K)$ is such that $i^\ast(e)=e_0$. The
topological group $\Gamma_0$ is isomorphic to the product space 
$G\times K_0$ with multiplication
$$
(g_1,a)(g_2,b)=(g_1g_2,m_{0}(g_1,g_2)+g_2 a+b).
$$
The continuous 2-cocycles $i^\ast (m_K)(g_1,g_2):=i(m_K
(g_1,g_2))$
and  $m_{0}$ are in the same class in $H^2_c(G,K_0)$. Then
$\Gamma_0$ is isomorphic to the extension $G\times K_0$ which has
multiplication
$$(g_1,k_1)(g_2,k_2)=(g_1g_2,i^\ast (m_{K})(g_1,g_2)+ g_2 k_1+k_2),$$
where $g_1,g_2\in G$ and $k_1,k_2\in K_0$. Let
$\widetilde{\Gamma}=G\times K$ be the group with multiplication
$$(g_1,k_1)(g_2,k_2)=(g_1g_2,m_{K}(g_1,g_2)+ g_2 k_1+k_2)$$
and
$\theta :G\times K\rightarrow G\times K_0$ be the map defined by
 $\theta(g,k)=(g,i(k))$. It is easy to see that $\theta$ is an injective homomorphism from
$\widetilde{\Gamma}$ to $\Gamma_0$. Let
$\Gamma:=\theta(\widetilde{\Gamma})$, then $\Gamma$ is a subgroup
of $\Gamma_0$ which is an extension of $K$ by $G$. Let $(1,i(K))$ be the subgroup of $\Gamma$ of elements $\{(1,i(k))\,:\, k\in K\}$. Analogously we use the notation $(G,0)$ and $(1,K_0)$. It is simple to
verify that  the intersection $\Gamma\cap (1,K_0)=(1,i(K))$ and $\Gamma
(1,K_0)=\Gamma_0$. Finally we show that $\Gamma$ is closed in  $\Gamma_0$. 
By hypothesis  $i(K)$ is closed in $K_0$ and  so $i(K)$ is the complement in $K_0$ of an open set $A$.  The set $A_0:=G \times A$ is an open set in $\Gamma_0$, since $\Gamma_0$ is equipped of
the product topology. Hence $\Gamma$, which is  equal to the complement in $\Gamma_0$ of $A_0$, is closed in $\Gamma_0$. 
\cvd

\def\FF{\mathbb{F}}

\section{Modules and cohomology groups for $\Sym(\Omega)$}\label{CG}
Henceforth, $\mathbb{F}$ will denote an arbitrary finite field  and $\mathbb{F}_2$ is  
the field with 2 elements, which (somewhat abusively) we identify with $\Z_2$. Suppose $G$ is a permutation group on a set $W$. Then $G$ acts on  $\FF^W$ the set of functions from $W$ to $\FF$ and if we give $\FF^W$ the product topology (with $\FF$ being discrete) then it is a profinite continuous $G$-module. We can also consider the (discrete) $G$-module $\FF W$, the $\FF$-vector space with basis $W$. There is a natural pairing $\FF^W \times \FF W \to \FF$ given by $(f, \sum_w a_w w) \mapsto \sum_w a_w f(w)$. By a standard application of Pontryagin duality, the  closed $\FF G$-submodules of $\FF^W$ are of the form $X^0$ for $\FF G$-submodules $X$ of $\FF W$, where $X^0$ denotes the annihilator of $X$ with respect to this pairing. Moreover if $X \leq Y$ are $\FF G$-submodules of $\FF W$ then $Y^0 \leq X^0$ and $X^0/Y^0$ is isomorphic to the Pontryagin dual $S^*$ of $S=Y/X$.

From now on, we take $\Omega$ to be an infinite set and let $G =\Sym(\Om)$. If $k$ is a positive integer 
then $\OmK$ denotes the set of subsets of size $k$ from $\Om$.  We let $W = \OmK$. We now collect together various results about closed submodules of $\mathbb{F}^W$ and associated cohomology groups.

As  is shown in \cite{Gr} (see also \cite{Jam}), the submodule structure of $\mathbb{F}\OmK$ is completely determined by the  maps $
\beta_{k,j}:\mathbb{F}[\Om]^k \rightarrow   \mathbb{F}[\Om]^j
$ for $0 \leq j \leq k$, where
 $\beta_{k,j}(w)=\sum\{w': w' \in
[w]^j\}$ for $w\in [\Om]^k$ and then extended linearly. Every proper submodule of $\mathbb{F}\OmK$ is an  intersection of  kernels of these maps. The \textit{Specht module} $S^{k}$ can be characterized as $\bigcap_{i=0}^{k-1}\Ker\beta_{k,i}$  and is irreducible. Furthermore, $\mathbb{F}\OmK$ has a finite $\FF G$-composition series in which the composition factors are isomorphic to $S^l$ for $l=0,\dots,k$, each appearing exactly once.

The dual of $\beta_{k,j}$ is the map $\alpha_{j,k} : \mathbb{F}^{[\Om]^j} \rightarrow   \mathbb{F}^{[\Om]^k}$ given by $\alpha_{j,k}(f)(w) = \sum_{v \in [w]^j} f(v)$, for $f \in \mathbb{F}^{[\Om]^j}$ and $w \in [\Omega]^k$. The image of $\alpha_{j,k}$ is the annihilator of $\Ker \beta_{k,j}$. Dualizing Gray's results, every closed $\FF G$-submodule of $\mathbb{F}^{[\Om]^k}$ is a sum of images of various $\alpha_{j,k}$ (for $j \leq k$) and  the (topological) $\FF G$-composition factors of $\FF^\OmK$ are the duals $(S^j)^*$ of the Specht modules (for $j = 0, 1, \ldots, k-1$). 

The paper \cite{Gr} also gives an algorithm for computing the lattice of submodules of $\mathbb{F}[\Om]^k$: what is needed is to  check the divisibility by the characteristic of $\FF$  of a finite number of binomial coefficients. Using this method and applying duality we record the  lattices of the closed submodules of $\Fi_2^{[\Om]^k}$ for $k=1,2$.

The only proper closed $G$-submodule of $\Fi_2^\Om$ is the submodule of the constant functions $\Fi_2 = \im\alpha_{0,1}\cong (S^0)^*$ and $\Fi_2^\Om/\Fi_2 \cong (S^1)^\ast$. 

For $\Fi_2^{[\Om]^2}$ we have proper closed submodules $(S^2)^0 = \im\alpha_{1,2} \oplus \Fi_2$, $\im\alpha_{1,2}$, $\Fi_2 = \im\alpha_{0,2}$ and $0$. We have $\Fi_2^{[\Om]^2}/(S^2)^0 \cong (S^2)^*$ and $\im\alpha_{1,2} \cong (S^1)^*$.







We now give some calculations of continuous cohomology groups of the infinite symmetric group with coefficients in certain sections of $ \Fi_2^{\OmK}$ (for $k=1,2$) that will be useful in the sequel. 

Let $S$ be a set and $E$ be a group of permutations acting transitively on $S$ with the topology given in the introduction. Let $H:=E_w$ be the stabilizer in $E$ of a point $w\in S$. If $F$ is a finite abelian group considered as a trivial $H$-module and $M$ is the profinite $E$-module $F^W$   then the continous analogue of Shapiro's lemma (\cite{EHew}, Proposition 2.2)
implies that$H^n_c(E,M)=H^n_c(H,F)$  for every $n\in \mathbb{N}$.
Indeed, the $E$-module $M$ is a coinduced module from a finite module for the stabilizer of a point. 

\begin{lemma}\label{31}
Suppose $\Om$ is infinite and let $G=\Sym(\Om)$.  Consider
$\mathbb{F}$ as a trivial $G$-module. Then \begin{enumerate}
\item[(i)]\label{vanificazione}$H^n_c(G,\Fi)=0$
for all $n>0$.
\item[(ii)]\label{split ext}
 $H^n_c(G,\Fi^\Om)=0$ for every
$n>0$.
\end{enumerate}
\end{lemma}
\emph{Proof.} (i) This is by \cite{EHew} (Theorem 3.5). 

(ii) Let $a \in \Om$. By Shapiro's Lemma
$H^n_c(G,\Fi^\Om)=H^n_c(G_{a},\mathbb{F})$ and
as $G_a$ is isomorphic to $\Sym(\Om)$ the result then follows from (i).\cvd\\

We note without proof that Theorem 3.5 of \cite{EHew} has the following generalization.  Its main consequence here is  Corollary \ref{cong e coho} below. Although this will not be used in the sequel, it seems to us to be worth recording here.  For definitions and  notation we refer the reader to \cite{EHew}: what is important for applications is to that $G = \Sym(\Om)$ acting on $\Delta = \Om$ has a smooth strong type, and that the topology on the symmetric group is the same whatever $\OmK$ it is considered as acting on.

\begin{lemma}\label{inflaction-restriction}
Suppose $G$ is a closed permutation group on $\Delta$ with a smooth strong type $p$. Suppose $Q$ is a finite group and $F$ a finite abelian group (regarded as a trivial $G$-module and $Q$-module). Then, for $n\geq 0$
$$
H^n_c(G\times Q,F)=H^n(Q,F).
$$\cvd
\end{lemma}

\begin{cor}\label{cong e coho}
Let $G=$Sym$(\Om)$. Then for every $k\in \N$ and $n \geq 1$
$$
H^n_c(G,\Fi^{[\Om ]^k})=H^n(\textrm{Sym}_k,\Fi)
$$
 \end{cor}
 where $\Sym_k$ is the finite symmetric group of degree $k$.
 \emph{Proof.}
By Shapiro's lemma 
$$
H^n_c(G,\Fi_p^{[\Om ]^k})=H^n_c(\Sym(\Om\setminus \{1, \dots, k\})\times \Sym_k,\Fi_p)
$$  
Applying  Lemma \ref{inflaction-restriction} to $ \Sym(\Om\setminus \{1, \dots, k\})\times \Sym_k$ we get the result.\cvd


%


\medskip

The following is Corollary 3.11 of \cite{EG}.

\begin{prop}\label{gray coho}
Let $G=\Sym(\Om)$ and $K_0=\Fi_2^{[\Om]^k}$.Then
\begin{displaymath}
H^1_c(G,(S^l)^\ast)=\left\{ \begin{array}{ll}
0 & \textrm{if $l\neq 2$}\\
\mathbb{Z}_2 & \textrm{if $l=2$.}
\end{array} \right.
\end{displaymath}
Moreover, If $K$ is a closed $G$-submodule of $K_0$, then
$H^1_c(G,K_0/K)$ is trivial if $(S^2)^\ast$ is not a composition
factor of $K_0/K$, otherwise it has order 2.\cvd
\end{prop}
\begin{lemma}\label{coho s1}
Let $G=Sym(\Om)$ and $\Fi = \Fi_2$. Then
$$
H^n_c(G,(S^1)^\ast)=0
$$
for $n\geq 1$.
\end{lemma}
\emph{Proof.} Take the short exact sequence
$$
0\rightarrow \mathbb{F}_2 \rightarrow
\mathbb{F}_2^{\Om}\rightarrow (S^1)^\ast \rightarrow 0.
$$
We have the long exact sequence
$$
\cdots \rightarrow H^n_c(G,\mathbb{F}_2) \rightarrow
H^n_c(G,\mathbb{F}_2^\Om)\rightarrow
H^n_c(G,(S^1)^\ast)\rightarrow H^{n+1}_c(G,\mathbb{F}_2)
$$
and by Lemma \ref{split ext} $H^n_c(G,\mathbb{F}_2^\Om)= 0$ for every
$n \geq 1$. By Lemma \ref{vanificazione}
$H^{n+1}_c(G,\mathbb{F}_2)=0$ and hence 
$
 H^n_c(G,(S^1)^\ast) =0.
$
\cvd

\medskip

Only the first part of the following will be needed in the sequel, but the second statement provides a complement to the above.

\begin{lemma}\label{coho s2}
Let $S^2$ be the Specht module of $\Fi_2[\Om]^2$ and $(S^2)^0 = \im\alpha_{1,2} \leq \Fi_2^{[\Om]^2}$. 
Then
$
H^n_c(G, (S^2)^0) = \{0\}$  and $$ 
H^n_c(G,(S^2)^\ast) = H^n_c(G,\Fi_2^{[\Om]^2})= H^n(\Fi_2, \Fi_2)
$$
for $n\geq 1$.
\end{lemma}
\emph{Proof.} Using the description of the submodules of $\Fi_2^{[\Om]^2}$ we have a short exact sequence
\[ 0 \to \Fi_2\oplus (S^1)^* \to \Fi_2^{[\Om]^2} \to (S^2)^* \to 0\]
and $(S^2)^0 \cong \Fi_2\oplus (S^1)^*$. Thus 
$$H^n_c(G,(S^2)^0) \cong  H^n_c(G, \Fi_2) \oplus H^n(G, (S^1)^*) = \{0\}$$
using Lemmas \ref{31} and \ref{coho s1}. 

Using this and the portion of the long exact sequence
$$
\rightarrow H^n_c(G,(S^2)^0)\rightarrow H^n_c(G,\Fi_2^{[\Om]^2})\rightarrow
H^n_c(G, (S^2)^\ast)\rightarrow H^{n+1}_c(G,(S^2)^0) \to
$$ we get that  $H^n_c(G,(S^2)^\ast) = H^n_c(G,\Fi_2^{[\Om]^2})$. The remaining statement then follows from Corollary \ref{cong e coho}. \cvd

\section{Full subgroups of $\Gamma_k$ }

Throughout, $\Om$ will be an infinite set, $k \geq 2$ and $G = \Sym(\Om)$. 
As in Definition \ref{ckdef}, we  consider  the free finite cover of  $\langle [\Om]^k;\Sym(\Om)\rangle$  with  binding group $\Z_2$, fibre group $\mathbb{Z}_4$ acting regularly, and canonical homomorphism given by the sign function. Existence and uniqueness of this comes from Lemma 2.1.2 of \cite{EMI} and we denote it by $\langle C_k, \Gamma_k \rangle$. We shall classify all minimally full subgroups of the $\Gamma_k$: in the case $k=2$ we show that $\Gamma_2$ is itself minimal (Theorem \ref{application1}); for $k \geq 3$ the result is that all minimally full subgroups arise from a homomorphic image of $\Gamma_2$  inside $\Gamma_k$ (Theorem \ref{47}).

\subsection{The case $k = 2$}

The following answers Problem 8.8 posed in \cite{EMI}.

\begin{theorem}\label{application1}
Suppose $\Gamma \leq \Gamma_2$ is full. Then $\Gamma  = \Gamma_2$.
\end{theorem}
\emph{Proof.}  The   kernel of $\langle C_2;\Gamma_2 \rangle$ is $K_0 = \Fi_2^{[\Om]^2}$. Suppose that $\Gamma < \Gamma_2$ is full. Then $K = \Gamma\cap K_0$ is a closed $G$-submodule of $K_0$ and so by the description of the closed submodules of $K_0$ already given, $K$ is  contained in $K_1 = (S^2)^0$. Thus, by considering $K_1\Gamma$,  we may assume $K = K_1$.
The extension $0 \to K_0 \to \Gamma_2 \to G \to 1$ is non-split (essentially because $\Z_4$ does not split over $\Z_2$: cf. 2.1.5 of \cite{EMI}). Thus the $2$-cocycle class $e \in H^2_c(G,K_0)$ to which it gives rise is non-zero. Now, by Lemma \ref{coho s2}, $H^2_c(G, K_1) = \{0\}$. Thus by Theorem \ref{main theorem}, there can be no full subgroup $\Gamma$ of $\Gamma_2$ with $\Gamma \cap K_0 = K_1$. \cvd



\subsection{Submodule structure of $\Fi_2^{[\Om]^k}$}
\def\l{\ell}
\def\notdiv{\not\vert}
\def\im{{\rm im\,}}
Recall that for $k \leq \l$ we have defined  $\alpha_{k,\l} : \Fi_2^{[\Om]^k} \to \Fi_2^{[\Om]^\l}$ by $\alpha_{k,\l}(f)(w) = \sum_{y \in [w]^k} f(y).$ This is the Pontryagin dual of $\beta_{\l, k}$ and $\im \alpha_{k,\l} = (\ker\beta_{\l,k})^0$. Note that if $j \leq k \leq \l$ then 
\[ \alpha_{k,\l}\circ\alpha_{j,k} = \binom{\l-j}{k-j}\alpha_{j, \l}.\]

By the dual of Gray's results in \cite{Gr}, any closed $G$-submodule $K$ of $\Fi_2^{[\Om]^\l}$ is a sum of certain submodules $\im \alpha_{k,\l}$ (with $k \leq \l$) and $(S^2)^*$ is \textit{not} a composition factor of $\im \alpha_{k,\l}$ iff $2 \vert \binom{\l-2}{k-2}$. 

\begin{lemma} \label{41} The submodule $M = \sum\{\im \alpha_{k,\l} : k \leq \l \mbox{ and } 2 \vert \binom{\l-2}{k-2}\}$ is the largest closed submodule of $\Fi_2^{[\Om]^\l}$ which does not have $(S^2)^*$ as a composition factor. 
\end{lemma}

\noindent\textit{Proof.\/} We have $\alpha_{k,\l}\circ\alpha_{2,k} = \binom{\l-2}{k-2}\alpha_{2, \l}$. Thus if $2\notdiv\binom{\l-2}{k-2}$ then $\im \alpha_{k,\l} \geq \im \alpha_{2,\l}$, and the latter has $(S^2)^*$ as a composition factor.

On the other hand, if $2\vert\binom{\l-2}{k-2}$ then $\im \alpha_{2,k} \leq \ker\alpha_{k,l}$. By Gray's results, $\im \alpha_{2,k}$ has $(S^2)^*$ as a composition factor, and so $\im\alpha_{k,\l}$ does \textit{not} have $(S^2)^*$ as a composition factor. Thus $M$ does not have $(S^2)^*$ as a composition factor.

Conversely suppose $K$ is a closed submodule which does not have $(S^2)^*$ as a composition factor. We can write $K = \sum_{k \in I} \im\alpha_{k,\l}$ for some set $I$ and by the first paragraph, if $k \in I$ then $2 \vert \binom{\l-2}{k-2}$. Thus $K \leq M$. \cvd

\begin{lemma} \label{42} If a closed submodule $K$ of $\Fi_2^{[\Om]^\l}$ has $(S^2)^*$ as a composition factor, then $K \geq \im\alpha_{2,\l}$. 
\end{lemma}

\noindent\textit{Proof.\/} As in the proof of the above, $K \geq \im \alpha_{k,\l}$ for some $k \leq \l$ with $2 \notdiv \binom{\l-2}{k-2}$. But then $\alpha_{2,\l}$ factors through $\alpha_{k,\l}$, as above, so $K \geq \im\alpha_{2,\l}$. \cvd

\subsection{Explicit construction of the free covers $C_k$}

\def\Z{\mathbb{Z}}
Let $\Om$ be any infinite set. It will be convenient to have an explicit description of the free finite covers $\pi_k : C_k \to [\Om]^k$. We fix some ordering $\leq$ on $\Omega$ and let $\Z_4 = \{0,1,2,3\}$ be the additive group of integers modulo 4. Inside the latter, we identify $\Fi_2$ with $\{0, 2\}$. 

\def\eps{\epsilon}

\begin{defn} \rm For $k \geq 2$, define functions $\eps_k : \Sym\Omega \times [\Omega]^k \to \Z_4$ as follows.

\begin{enumerate} 
\item ($k = 2$) If $g \in \Sym(\Om)$ and $w = \{w_1, w_2\} \in [\Om]^2$ with $w_1 < w_2$, then  
\[\eps_2(g,w) = \left\{ \begin{array}{ll} 
									0 & \mbox{ if $gw_1 < gw_2$} \\
									1 &  \mbox{ otherwise}
									\end{array}
									\right. .\]
									
\item Suppose $k \geq 3$. Define $\eps_k : \Sym\Omega \times [\Omega]^k \to \Z_4$ by $$\eps_k(g,w) = \sum_{y \in [w]^2} \eps_2(g,y).$$

\end{enumerate}

\end{defn}

\def\ldot{.}

\noindent Now for $k \geq 2$ let 
\[\Gamma_k = \Fi_2^{[\Om]^k} \times \Sym(\Om)\]
and
\[C_k = \Z_4 \times [\Om]^k.\]
Define $\mu_k : \Gamma_k \times C_k \to C_k$ by 

\[\mu_k((f,g), (a,w)) = (a+f(g(w)) + \eps_k(g,w), gw).\]
Thus any $(f,g) \in \Gamma_k$ defines a map $\mu_k((f,g), \ldot): C_k \to C_k$, which can easily be seen to be a bijection. Computing the composition of two such maps, we obtain
\[\mu_k((\l, h),\ldot) \circ \mu_k((f,g), \ldot) = \mu_k((\l+{}^hf + c_k(h,g), hg), \ldot)\]
where 
\[c_k(h,g)(w) = \eps_k(g,w)+\eps_k(h,gw)-\eps_k(hg, w)\]
and $({}^hf)(w) = f(h^{-1}w)$. Note that $c_k(h,g) \in \Fi_2^{[\Om]^k}$. 

Thus if we define a product on $\Gamma_k$ by 
\[(\l, h)(f,g) = (\l+{}^hf+ c_k(h,g), hg)\]
then $\Gamma_k$ is a group and $\mu_k$ is a faithful action of $\Gamma_k$ on $C_k$, so we may regard $\Gamma_k$ as a subgroup of $\Sym(C_k)$. 

Define $\pi_k : C_k \to [\Om]^k$ by $\pi_k((a,w)) = w$. It is clear that $\Gamma_k$ preserves the fibres of $\pi_k$ and the kernel of the action of $\Gamma_k$ on the fibres is $\Fi_2^{[\Om]^k} \times 1 $. We shall show below that $\Gamma_k$ is closed in $\Sym(C_k)$. Thus $\langle C_k; \Gamma_k\rangle$ is a permutation structure and $\pi_k : C_k \to [\Om]^k$ makes it a finite cover of $\langle [\Om]^k, \Sym(\Om)\rangle$. It is easily checked (by direct calculation of the stabilizer of $w \in [\Om]^k$ in $\Gamma_k$) that the fibre group is $\Z_4$ and the binding group is $\Fi_2$ Assuming that $\Gamma_k$ is closed in $\Sym(C_k)$, we therefore have:

\begin{prop} The permutation structure $\langle C_k; \Gamma_k\rangle$ is the free finite cover of $\langle [\Om]^k; \Sym(\Om)\rangle$ with these data.
\end{prop}

So it remains to prove that $\Gamma_k$ is closed in $\Sym(C_k)$. Note that if $K$ is a closed normal subgroup of a topological group $\Pi$ and $H$ is a subgroup of $\Pi$ such that $K \cap H$ is closed, then $H$ is closed in $\Pi$ if $\overline{H}\cap K \subseteq H$ (where $\overline{H}$ is the closure of $H$ in $\Pi$). We want to apply this where $H = \Gamma_k$, $\Pi$ is the permutations of $C_k$ which preserve all fibres of $\pi_k$, so a closed subgroup of $\Sym(C_k)$, and $K = \Z_4^{[\Om]^k}$. We need to verify that any element of $K$ which lies in the closure of  $H$ is in $\Fi_2^{[\Om]^k}$. But this follows from the fact that if $A$ is a finite subset of $\Omega$ and $(f,g) \in \Gamma_k$ fixes every element of $A$, then $f \vert_{ [A]^k} \in \Fi_2^{[A]^k}$. 

\subsection{Lifting the $\alpha_{k,\l}$}

\begin{prop} For $k \geq 2$, define $\gamma_{2,k} : \Gamma_2 \to \Gamma_k$ by 
$\gamma_{2,k}(f,g) = (\hat{f}, g)$, where $\hat{f} = \alpha_{2,k}(f)$. Then $\gamma_{2,k}$ is a continuous group homomorphism which extends $\alpha_{2,k} : \Fi_2^{[\Om]^2} \to \Fi_2^{[\Om]^k}$, and $\im \gamma_{2,k}$ is closed in $\Gamma_k$. 

\end{prop}

\noindent\textit{Proof.\/} It is clear that $\gamma_{2,k}$ extends $\alpha_{2,k}$. To check that it is a group homomorphism we need to show that for $(\l,h), (f,g) \in \Gamma_2$ we have:
\[\alpha_{2,k}(\l+{}^hf + c_2((h,g)) = \hat{\l} + {}^h\hat{f} + c_k(h,g).\]
As $\alpha_{2,k}$ is a $G$-homomorphism this reduces to showing that $\alpha_{2,k}(c_2(g,h) )= c_k(g,h)$, and this follows trivially from the definition of $\eps_k(h,g)$ and $c_k(h,g)$. 

Continuity of $\gamma_{2,k}$ is routine and closedness of the image will be part of the next proposition. \cvd

\begin{prop} \label{46} Suppose $2 \leq k \leq \l$ and $2 \not\vert \binom{\l-2}{k-2}$. Then $\gamma_{k,\l} : \Gamma_k \to \Gamma_\l$ given by $\gamma_{k,\l}((f,g)) = (\hat{f}, g)$, where $\hat{f} = \alpha_{k,\l}(f)$, is a continuous group homomorphism which extends $\alpha_{k,\l}$ and whose image $\im\gamma_{k,\l}$ is closed in $\Gamma_\l$.

\end{prop}

\noindent\textit{Proof.\/} Similarly to the previous result, we need to show that $\alpha_{k,\l}(c_k(h,g)) = c_\l(h,g)$. Now, for $v \in [\Om]^\l$
\[c_\l(h,g)(v) = \sum_{y \in [v]^2}(\eps_2(g,y) + \eps_2(h, gy) - \eps_2(hg, y).\]
Also 
\[\alpha_{k,\l}(c_k(h,g)(v)) = \sum_{w \in [v]^k}(\eps_k(g,w) + \eps_k(h, gw) - \eps_k(hg, w))\]
and this is equal to $\binom{\l-2}{k-2} c_\l(g,h)(v)$. As $2 \not\vert \binom{\l-2}{k-2}$ this is equal to $c_\l(g,h)$, because we are working in $\Fi_2$. 

Closedness of the image is as at the end of the previous subsection, using the fact that if $g$ fixes every element of some finite subset $A$ of $\Om$ and $w \in [A]^k$ then $\eps_k(g,w) = 0$. \cvd

\subsection{Minimally full subgroups of $\Gamma_k$}

The main result is:

\begin{theorem}\label{47}  Suppose $k \geq 3$. Then the conjugates of $\im \gamma_{2,k}$ are the minimally full subgroups of $\Gamma_k$.
\end{theorem}

The proof is a series of lemmas. Let $K_0 \leq \Gamma_k$ be the kernel $\Fi_2^{[\Om]^k}$ of the cover $C_k$. If $\Gamma \leq \Gamma_k$ is full we refer to $K_0 \cap \Gamma$ as the \textit{kernel} of $\Gamma$.

\begin{lemma} The closed subgroup $\im \gamma_{2,k}$ is a minimally full subgroup of $\Gamma_k$ with kernel $\im \alpha_{2,k}$.

\end{lemma}

\noindent\textit{Proof.\/} It remains to see that $\im \gamma_{2,k}$ is minimal. But by Theorem \ref{application1}, $\Gamma_2$ has no proper closed subgroup which projects onto $\Sym(\Om)$. The same is therefore true of any continuous homomorphic image of $\Gamma_2$. \cvd

\begin{lemma} \label{49} Suppose $K$ is closed and $\im\alpha_{2,k} \leq K \leq K_0$. Then there is exactly one $\Gamma_k$-conjugacy class of full subgroups $H$ of $\Gamma_k$ with $H \cap K_0 = K$. 

\end{lemma}

\noindent\textit{Proof.\/} To see that there is some such subgroup, note that $H = (\im \gamma_{2,k})K$ is closed and satisfies $H \cap K_0 = K$. To see that any two such closed subgroups are conjugate, it suffices to show that $H^1_c(G, K_0/K) = \{0\}$. This follows from Proposition \ref{gray coho} as $(S^2)^*$ is a composition factor of $K$, and therefore not a composition factor of $K_0/K$. \cvd

\begin{lemma}\label{410}  Suppose $H$ is a minimally full subgroup of $\Gamma_k$ with kernel $K$ and $K$ has $(S^2)^*$ as a composition factor. Then $K = \im \alpha_{2,k}$ and $H$ is a $\Gamma_k$-conjugate of $\im \gamma_{2,k}$. 
\end{lemma}

\noindent\textit{Proof.\/} By Lemma \ref{42},  $K \geq \im\alpha_{2,k}$. By Lemma \ref{49}, $H$ contains a $\Gamma_k$-conjugate of $\im\gamma_{2,k}$. Therefore minimality implies that $H$ is a $\Gamma_k$-conjugate of $\im\gamma_{2,k}$. \cvd

\medskip

\noindent\textit{Proof of Theorem \ref{47}:\/} By Lemma \ref{410} it is enough to prove that if $K \leq K_0$ does not have $(S^2)^*$ as a compositon factor, then $K$ is not the kernel of a full subgroup of $\Gamma_k$. So suppose for a contradiction that  $H$ is a full subgroup of $\Gamma_k$ whose kernel $K = H \cap K_0$ does not have $(S^2)^*$ as a composition factor. By Lemma \ref{41}
\[ K \leq M = \sum\{\im \alpha_{j,k} : j \leq k \mbox{ and } 2 \vert \binom{k-2}{j-2}\}\]
and so there is a full subgroup $HM \leq \Gamma_k$ with kernel $M$. Thus, for our contradiction we may assume that $K = M$.

Let $\eta : M \to K_0$ be inclusion and $e \in H^2_c(G, K_0)$ the $2$-cocycle class arising from the extension $0 \to K_0 \to \Gamma_k \to G \to 1$. Note that the construction of $\Gamma_k$ gives an explicit continuous $2$-cocycle in $e$, namely $c_k$. In the notation of Theorem \ref{main theorem}, the existence of $H$ implies that in the portion of the long exact sequence

\[ \cdots \to H^2_c(G, M) \stackrel{\eta^*}{\to} H^2_c(G,K_0) \stackrel{\zeta}{\to} H^2_c(G, K_0/M) \to \cdots\]
the class $e$ is in the image of $\eta^*$, therefore it is in $\ker(\zeta)$. The following lemma will allow us to embed $K_0/M$ in some $\Fi_2^{[\Om]^\l}$.

\begin{lemma}\label{411} \begin{enumerate} \item There is $\l > k$ such that
\begin{enumerate}\item[(i)] if $j \leq k$ and $2 \vert \binom{k-2}{j-2}$ then $2\vert \binom{\l -j}{k-j}$; 
\item[(ii)] $2\not\vert\binom{\l-2}{k-2}$. 
\end{enumerate}
\item For $\l$ as in (1), $\ker\alpha_{k,\l} = M$. 
\end{enumerate}

\end{lemma}

Indeed, suppose $\l$ is as in the lemma and so $\alpha_{k,\l} : \Fi_2^{[\Om]^k} \to \Fi_2^{[\Om]^\l}$ has kernel $M$. By condition (ii) on $\l$ and Proposition \ref{46}, $\alpha_{k,\l}$ extends to a homomorphism $\gamma_{k,\l} : \Gamma_k \to \Gamma_\l$ which also has kernel $M$. Now, $\im \gamma_{k,\l}$ is a covering expansion of $\Gamma_\l$ with kernel isomorphic to $K_0/M$, and it is easily checked that the fibre group in this expansion is $\Z_4$ and the binding group $\Fi_2$. In particular, the fibre group does not split over the binding group and it follows (cf. \cite{EMI}, Proof of Lemma 2.1.5) that the extension
\[0 \to K_0/M \to \im\gamma_{k, \l} \to G \to 1\]
is non-split. In particular the continuous $2$-cocycle class in $H^2_c(G, K_0/M)$ corresponding to this extension is non-zero. but this class is $\zeta(e)$, so we have a contradiction. \cvd

\medskip

It remains to prove Lemma \ref{411}. To do this we recall some information from (\cite{Jam}, pp. 87--88) about binomial coefficients.

\begin{lemma}\label{413} Suppose $b \leq a \in \N$ are written $2$-adically as 
\[a = \sum_{i \leq r} a_i 2^i \mbox{ and } b = \sum_{i \leq r} b_i2^i\]
with $a_i, b_i \in \{0,1\}$. Then 
\[ 2 \vert \binom{a}{b} \Leftrightarrow a_i < b_i \mbox{ for some } i \leq r. \,\,\cvd\]

\end{lemma}

Now we prove Lemma \ref{411}. For part (1), write $a = k-2$ and $e = \l - k$. We want to choose $e$ so that 
\begin{enumerate}
\item[(i)] for all $b \geq 1$, if $2 \vert \binom{a}{b}$ then $2 \vert \binom{e+b}{e}$;
\item[(ii)] $2 \not\vert \binom{e+a}{e}$.
\end{enumerate}
\def\ba{\bar{a}}
Suppose $a = a_r2^r+\ldots+a_02^0$ (with $a_i \in \{0,1\}$ and $a_r \neq 0$). So in binary notation, $a = a_ra_{r-1}\ldots a_0$. Let $e = 2^{r+1}+\ba_r2^r+\ba_{r-1}2^{r-1}+\ldots + \ba_0$ where $\ba_i = 1 \Leftrightarrow a_i = 0$. So in  binary notation, $e = 1\ba_r\ba_{r-1} \ldots \ba_0$.  We claim that this works.

For (ii) note that $a+e = 11\ldots1$ (in binary notation), so by Lemma \ref{413}, $2 \not\vert \binom{a+e}{e}$. For (i) suppose $a > b = b_rb_{r-1}\ldots b_0$ (in binary) and $2 \vert\binom{a}{b}$. Then $b_i > a_i$ for some $i \leq r$, which we may take to be a small as possible. Note that $b_i = 1$ and $a_i = 0$. We now compare the binary digits $(e+b)_i$ and $b_i$ of $e+b$ and $b$ respectively. Note that $e_i = \ba_i = 1$. Furthermore, if $j < i$ and $e_j = b_j = 1$, then $a_j = 0 < b_j$, contradicting the minimality of $i$. It follows that $(e+b)_i = 0 < 1 = e_i$, and so $2 \vert \binom{e+b}{e}$, as required.

To prove part (2) of the lemma, suppose $\l$ satisfies the conditions in part (1). Recall that for $2 \leq j \leq k \leq \l$ we have
\[\alpha_{k,\l}\circ\alpha_{j,k} = \binom{\l-j}{k-j}\alpha_{j,\l}.\]
If $2 \vert \binom{k-2}{j-2}$ then $2 \vert \binom{\l-j}{k-j}$, and therefore $\im \alpha_{j,k} \leq \ker\alpha_{k,\l}$. So $M \leq \ker\alpha_{k,\l}$.

On the other hand, as $2 \not\vert \binom{\l-2}{k-2}$ we have $\alpha_{k,\l}\circ \alpha_{2,k} = \alpha_{2,\l} \neq 0$. So $\im \alpha_{2,k} \not\leq \ker\alpha_{k,\l}$ and therefore by Lemma \ref{42}, $\ker\alpha_{k,\l}$ does not have $(S^2)^*$ as a composition factor. It then follows from Lemma \ref{41} that $\ker\alpha_{k,\l} = M$, as required. \cvd


\begin{thebibliography}{LGO \,}
\normalsize \setlength{\parskip}{6pt}




\bibitem{AZ2} G. Ahlbrandt and M. Ziegler,
\emph{What's so special about $(\mathbb{Z}/4\mathbb{Z})^\omega$?},
Archive for Mathematical Logic \textbf{31} (1991), 115--132.


\bibitem{Ca}  P.J. Cameron,  \emph{Oligomorphic Permutation Groups}, London Mathematical Society Lecture Notes Series 152, Cambridge
University Press, Cambridge (1990).

\bibitem{CKK} J. Cossey, O. H. Kegel and L. G. Kov\'{a}cs, \emph{Maximal Frattini extensions}, Arch. Math. (Basel) 35 (1980), 210--217.

\bibitem{EG} D.M. Evans and D.G.D. Gray, \emph{Kernels and cohomology groups for some finite covers},
Lecture Notes in Logic, Logic Colloquium 96, (Eds. J. M.
Larrazzabal, D. Lascar and G. Mints, Springer, 1998.

\bibitem{DEJAlg} D. M. Evans, \emph{Computation of first cohomology groups of finite covers}, J. Algebra 193 (1997), 214--238.

\bibitem{EMI}   D.M. Evans, A. A. Ivanov and D. Macpherson, \emph{Finite covers}, in   \emph{Model Theory of Groups and Automorphism Groups}, London Mathematical Society Lecture Notes Series 244, Cambridge
University Press, Cambridge (1977), \mbox{1--72} .



\bibitem{EHew} D.M. Evans and P.R. Hewitt, \emph{Continuous cohomology of permutation groups on profinite modules}, Communications in Algebra 34 (2006), 1251-- 1264.

\bibitem{Gr} D.G.D. Gray, \emph{The structure of some permutation modules for the symmetric group of infinite degree},
Journal of Algebra, \textbf{193} (1997), \mbox{122--143.}


\bibitem{HPi} W.A. Hodges and A. Pillay, \emph{Cohomology of structures and some problems of Ahlbrandt and
Ziegler}, J. London Math. Soc. (2) \textbf{50} (1994), 1--16.


\bibitem{Jam} G.D. James, \emph{The Representation Theory of the Symmetric Groups}, Springer-Verlag, 1978.



\end{thebibliography}
\end{document}